\documentclass[a4paper,12pt]{article}
\usepackage[top=2.5cm,bottom=2.5cm,left=2.5cm,right=2.5cm]{geometry}
\usepackage{amssymb}
\usepackage{graphicx}
\usepackage{amsmath,amsthm,amssymb,lineno}
\usepackage{latexsym}
\usepackage{graphicx,booktabs,multirow}
\usepackage{appendix}
\usepackage{longtable}
\graphicspath{{figures/}}

\newtheorem{theorem}{Theorem}

\newtheorem{corollary}[theorem]{\rm\bfseries Corollary}

\begin{document}
%\linenumbers

\title{Tutte polynomial of a fractal scale-free lattice
\thanks{Project supported by Hunan Provincial Natural Science
Foundation of China (13JJ3053).}}
\author{Hanlin Chen,Yuanhua Liao, Hanyuan Deng\thanks{Corresponding author:
hydeng@hunnu.edu.cn}\\
{\footnotesize College of Mathematics and Computer Science,} \\
{\footnotesize Hunan Normal University, Changsha, Hunan 410081, P. R. China} \\
}

\date{}
\maketitle

\begin{abstract}
The Tutte polynomial of a graph, or equivalently the $q$-state Potts
model partition function, is a two-variable polynomial graph
invariant of considerable importance in both combinatorics and
statistical physics. The computation of this invariant for a  graph
is NP-hard in general. In this paper, based on their self-similar
structures, we recursively describe the Tutte polynomials of an
infinite family of scale-free lattices. Furthermore, we give some
exact analytical expressions of the Tutte polynomial for several
special points at $(X,Y)$-plane.

%The main aim of this paper is to compute the Tutte polynomial of  a fractal scale-free lattice(network),
%
% Moreover, we also got the Tutte polynomial of two other scale-free networks,
% (2,2)-flower and (1,3)-flower, which much similar, and we can quickly get the spanning tree of these networks.

\noindent
%\textbf{AMS classification}: 05C30, 05C69\\
{\bf Keywords}: Potts model; Tutte polynomial; Spanning tree;
Asymptotic growth constant; Self-similar.
\end{abstract}

\baselineskip=0.30in

\section{Introduction}\label{intro}

The Tutte polynomial $T(G;x,y)$ of a graph $G$, due to W. T. Tutte
\cite{Tutte54}, is a two-variable polynomial and is deeply connected
with many areas of both  physics and mathematics. It is defined as
$$T(G;x,y)=\sum\limits_{H\subseteq G}(x-1)^{r(G)-r(H)}(y-1)^{n(H)},$$
where the sum runs over all the spanning subgraphs $H$ of $G$,
$r(G)=|V(G)|-k(G)$, $n(G)=|E(G)|-|V(G)|+k(G)$ and $k(G)$ is the
number of components of $G$.  For a thorough survey on the Tutte
polynomial, we refer the reader to
\cite{Jaeger90,Welsh00JMP,Ellis11}. The importance of the Tutte
polynomial comes from the rich information it contains about the
underlying graph. It contains several other polynomial invariants,
such as the chromatic polynomial, the flow polynomial, the Jones
polynomial and the all terminal reliability polynomial as partial
evaluations, and various numerical invariants such as the number of
spanning trees as complete evaluations. Furthermore, the Tutte
polynomial has been widely studied in the field of statistical
physics where it appears as the partition function of the zero-field
Potts model \cite{Potts52,Wu92}. In fact, let $G$ be a graph with
$n$ vertices and $k(G)$ components and let $Z_{zero}$ denote the
partition function of the zero-field Potts model, then
$Z_{zero}(G;q,v)$ and $T(G;x,y)$ satisfy the following relation
\cite{Fortuin72,Beaudin10}:
\begin{equation*}
Z_{zero}(G;q,v)=q^{k(G)}v^{n-k(G)}T(G;(q+v)/v,v+1).
\end{equation*}

In both fields of combinatorics and statistical physics, Tutte
polynomials of many graphs (or lattices) have been computed by
different methods
\cite{Shrock00a,Shrock00b,Shrock01,Chang04,Chang08JSP,Shrock11JPA,Shrock12,Salas09JSP,Alvarez12}.
Recently, on the basis of the subgraph expansion definition of the
Tutte polynomial, a very useful method for computing the Tutte
polynomial, called the subgraph-decomposition method, was developed
by Donno \emph{et al.} \cite{Donno13}. This technique is highly
suited for computing the Tutte polynomial of self-similar graphs,
and some applications of it can be found in
\cite{Liao13EPL,Liao13Physica,Gong14}.

The lattice (network) under consideration was introduced by Kaufman
\emph{et al.} \cite{Kaufman,Griffiths} and was further studied by
Zhang \emph{et al.} \cite{Zhang08JSM,Zhang11PRE} from the viewpoint
of complex networks. The aim of this paper is to compute the Tutte
polynomial of this deterministic scale-free network. For this
purpose, we first partition the set of spanning subgraphs of this
network into two disjoint subsets. In this way, we can express the
Tutte polynomial by two summands. Then, we study the contributions
of all kinds of possible combinations of this two subsets. Finally,
base on the self-similar structure of the network, we obtain a
recursive formula for each summand. In particular, as special cases
of the general Tutte polynomial, we get:

\begin{itemize}
  \item the recursive formula for computing the Tutte
polynomial $T(G;x,y)$;
  \item the number $\tau(G)$ of spanning trees and the asymptotic growth
constant
$\lim\limits_{n\rightarrow\infty}\frac{ln\tau(G_n)}{|v(G_n)|}$;
  \item the dimension
of the bicycle space;
  \item the number of acyclic root-connected orientations;
the number of indegree sequences of strongly connected orientations.
\end{itemize}

\section{Preliminaries}

In this section, we briefly discuss some necessary background that
will be used for our calculations. We use standard graph terminology
and the words ``network"  and ``graph" indistinctly.  Let $G$ be a
graph with vertex set $V(G)$ and edge set $E(G)$. The vertices $a$
and $b$ are the end-points of an edge $\{a,b\}$. An \emph{orientation}
 of graph $G$ is the digraph defined by the choice of a direction for
 every edge of $E(G)$. A \emph{directed cycle} of a digraph is a set of
  edges forming a cycle of the graph such that they are all directed
  accordingly with a direction for the cycle. A digraph is \emph{acyclic}
   if it has no directed cycle, and \emph{strongly connected} if for every
    pair of vertices there is a directed cycle containing them.
   A \emph{sink} for a digraph is a vertex with no outgoing edge.

 %%%%%%%%%%%%%%%%%%%%%%%%%%%%%%%%%%%%%%%%%%%%%%%%%%%%%%%%%%%%%%%%%%%%55
                          %Figure 1$
  %%%%%%%%%%%%%%%%%%%%%%%%%%%%%%%%%%%%%%%%%%%%%%%%%%%%%%%%%%%%%%%%%%%%55
  \begin{figure}[ht!]
\begin{center}
\includegraphics[width=10cm]{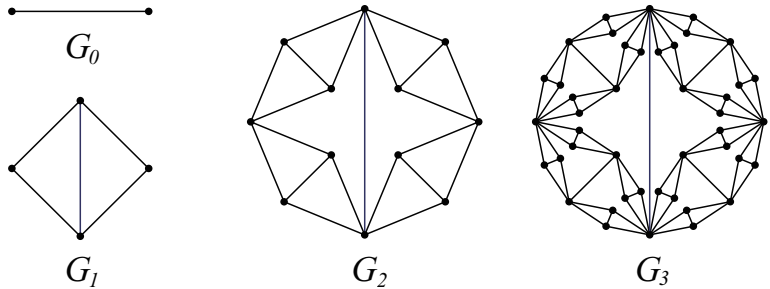}
\caption{The lattices $G_{0},G_{1},G_{2}$ and $G_{3}$.}
\end{center}
\end{figure}
 %%%%%%%%%%%%%%%%%%%%%%%%%%%%%%%%%%%%%%%%%%%%%%%%%%%%%%%%%%%%%%%%%%%%55
  %%%%%%%%%%%%%%%%%%%%%%%%%%%%%%%%%%%%%%%%%%%%%%%%%%%%%%%%%%%%%%%%%%%%55

The network, as shown in Figure 1, can be constructed as follows.

For $n=0$, $G_0$ is the complete graph $K_2$. For $n\geq 0$,
$G_{n+1}$ can be constructed from four copies of $G_{n}$ by merging
four groups of vertices and adding a new edge. Specifically, let
$X_n$ and $Y_n$, hereafter called special vertices of $G_n$, be the
leftmost and the rightmost vertex of $G_n$. $X_{n}$ and $X_{n}$ are
combined into the special vertex $X_{n+1}$ of $G_{n+1}$, $Y_{n}$ and
$Y_{n}$ are combined into the special vertex $Y_{n+1}$ of $G_{n+1}$,
and a new edge $e_n$ is added between two vertices combined by
$Y_{n}$ and $X_{n}$. The construction of $G_{n+1}$ is illustrated in
Figure 2.

%%%%%%%%%%%%%%%%%%%%%%%%%%%%%%%%%%%%%%%%%%%%%%%%%%%%%%%%%%%%%%%%%%%55
                          %Figure 2$
  %%%%%%%%%%%%%%%%%%%%%%%%%%%%%%%%%%%%%%%%%%%%%%%%%%%%%%%%%%%%%%%%%%%%55
  \begin{figure}[ht!]
\begin{center}
\includegraphics[width=10cm]{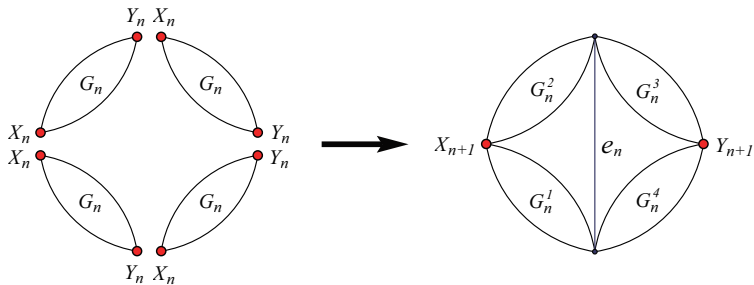}
\caption{The construction of $G_{n+1}$}
\end{center}
\end{figure}
 %%%%%%%%%%%%%%%%%%%%%%%%%%%%%%%%%%%%%%%%%%%%%%%%%%%%%%%%%%%%%%%%%%%%55
  %%%%%%%%%%%%%%%%%%%%%%%%%%%%%%%%%%%%%%%%%%%%%%%%%%%%%%%%%%%%%%%%%%%%55

\section{Recursive formulas}

In this section, we give the computational formulas of Tutte polynomial of the lattice $G_n$.
It is easy to obtain that the order and size of the lattice $G_n$ are, respectively,
$$|V(G_n)|=\frac{2\times4^n+4}{3},  ~~~~   |E(G_n)|=\frac{4^{n+1}-1}{3}.$$
And the average degree after $n$ iterations is $\langle $k$ \rangle_{n}=\frac{2|E(G_n)|}{|V(G_n)|}$,
which approaches 4 in the infinite $n$ limit.

To investigate the Tutte polynomial $T(G_n;x,y)$. First of all, we
partition the set of the spanning subgraph of $G_n$ into two
disjoint subsets:
\begin{itemize}
\item  $\mathcal{G}_{1,n}$ denotes the set of spanning subgraphs of $G_n$,
where two special vertices $X_n$ and $Y_n$ of $G_n$ belong to the
same component;
\item  $\mathcal{G}_{2,n}$ denotes the set of spanning subgraphs of $G_n$,
where two special vertices $X_n$ and $Y_n$ of $G_n$ do not belong to
the same component.
\end{itemize}

%%%%%%%%%%%%%%%%%%%%%%%%%%%%%%%%%%%%%%%%%%%%%%%%%%%%%%%%%%%%%%%%%%55
                          %Figure 3$
%%%%%%%%%%%%%%%%%%%%%%%%%%%%%%%%%%%%%%%%%%%%%%%%%%%%%%%%%%%%%%%%%%%%55
\begin{figure}[ht!]
\begin{center}
\includegraphics[width=7cm]{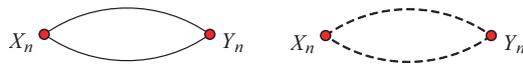}
\caption{The spanning subgraph of $G_{n}$: Type I (left), Type II
(right).}
\end{center}
\end{figure}
%%%%%%%%%%%%%%%%%%%%%%%%%%%%%%%%%%%%%%%%%%%%%%%%%%%%%%%%%%%%%%%%%%%%55
%%%%%%%%%%%%%%%%%%%%%%%%%%%%%%%%%%%%%%%%%%%%%%%%%%%%%%%%%%%%%%%%%%%55

Observe that, for each $n\geq 0$, we have the partition
$\mathcal{G}_{1,n}\cup \mathcal{G}_{2,n}$ of the set of spanning
subgraphs of $G_n$. Next, let $T_n(x,y)$ denote the Tutte polynomial
$T(G_n;x,y)$ of $G_n$. For every $n\geq1$, $T_{1,n}(x,y)$ and
$T_{2,n}(x,y)$ are the following polynomials:

\begin{itemize}
  \item $T_{1,n}(x,y)= \sum_{H \in \mathcal{G}_{1,n}}(x-1)^{r(G_n)-r(H)}(y-1)^{n(H)}$;
  \item $T_{2,n}(x,y)=\sum_{H \in \mathcal{G}_{2,n}}(x-1)^{r(G_n)-r(H)}(y-1)^{n(H)}$.
\end{itemize}
We have
\begin{equation*}
T_n(x,y)=T_{1,n}(x,y)+T_{2,n}(x,y).
\end{equation*}

In order to obtain $T_n(x,y)$, we need to find recursive formulas on
$T_{1,n}(x,y)$ and $T_{2,n}(x,y)$. For this purpose, we analyze the
relation between spanning subgraphs of $G_{n+1}$ and spanning
subgraph of $G_{n}$. Observe from Figure 2 that $G_{n+1}$ is
constructed from four copies of $G_{n}$ by merging some special
vertices and adding a new edge $e_n$. Thus, a spanning subgraphs of
$G_{n+1}$ is combined by the spanning subgraphs of the four copies
$G_{n}^{i} (i=1,2,3,4)$ of $G_{n}$ with $S$, where $S$ may be
$\{e_{n}\}$ or $\emptyset$. Indeed, a spanning subgraph $H$ of
$G_{n+1}$ is uniquely determined by the restriction of $H$ to the
four copies $G_{n}^{i}$ (denoted by $H_i(i=1,2,3,4)$, respectively)
and $S$, and vice versa. Therefore, the Tutte polynomial of
$G_{n+1}$ can be written as
$$T_{n+1}(x,y)=\sum_{H_i \subseteq G_{n}^{i};(\bigcup\limits_{i=1}^4 H_i) \bigcup S=H}(x-1)^{r(G_{n+1})-r(H)}(y-1)^{n(H)}.$$
where the sum runs over all spanning subgraphs $H_i$ of $G_{n}^{i}$
($i=1,2,3,4$) and $S$. Now, we need to know how $r(H)$ and $n(H)$
depend on $r(H_i)$ and $n(H_i)$, for $i=1,2,3,4$. Note that
$|V(G_{n+1})|=4|V(G_n)|-4$, $|E(H)|=\sum\limits_{i=1}^{4}|E(H_i)|+1$
if $S=\{e_n\}$ and $|E(H)|=\sum\limits_{i=1}^{4}|E(H_i)|$ if
$S=\emptyset$. So, there are two cases to be considered.

{\bf Case 1}.  $S=\{e_n\}$.

In this case, we consider the spanning subgraph $H$ of $G_{n+1}$,
which contains the new adding edge $\{e_n\}$, and
$|E(H)|=\sum\limits_{i=1}^{4}|E(H_i)|+1$.

{\bf Subcase 1}. If $k(H)=\sum_{i=1}^{4}k(H_{i})-3$, then
\begin{equation*}
r(H)=|V(H)|-k(H)=(4|V(G_{n})|-4)-(\sum_{i=1}^{4}k(H_{i})-3)=\sum_{i=1}^{4}r(H_{i})-1.
\end{equation*}
Moreover, we have
\begin{equation*}
n(H)=|E(H)|-r(H)=(\sum_{i=1}^{4}|E(H_i)|+1)-(\sum_{i=1}^{4}r(H_{i})-1)=\sum_{i=1}^{4}n(H_i)+2
\end{equation*}
and
\begin{equation*}
r(G_{n+1})-r(H)=(|V(G_{n+1})|-1)-(\sum_{i=1}^{4}r(H_i)-1)=\sum_{i=1}^{4}(r(G_{n})-r(H_i)).
\end{equation*}
Thus,
\begin{equation*}
(x-1)^{r(G_{n+1})-r(H)}(y-1)^{n(H)}=(y-1)^{2}\prod_{i=1}^{4}(x-1)^{r(G_{n})-r(H_{i})}(y-1)^{n(H_i)}.
\end{equation*}

{\bf Subcase 2}. If $k(H)=\sum\limits_{i=1}^{4}k(H_{i})-4$, then
\begin{equation*}
r(H)=|V(H)|-k(H)=(4|V(G_{n})|-4)-(\sum_{i=1}^{4}k(H_{i})-4)=\sum_{i=1}^{4}r(H_{i}).
\end{equation*}
Moreover, we have
\begin{equation*}
n(H)=|E(H)|-r(H)=(\sum_{i=1}^{4}|E(H_i)|+1)-(\sum_{i=1}^{4}r(H_{i}))=\sum_{i=1}^{4}n(H_i)+1
\end{equation*}
and
\begin{equation*}
r(G_{n+1})-r(H)=(|V(G_{n+1})|-1)-(\sum_{i=1}^{4}r(H_i))=\sum_{i=1}^{4}(r(G_{n})-r(H_i))-1.
\end{equation*}
Hence,
\begin{equation*}
(x-1)^{r(G_{n+1})-r(H)}(y-1)^{n(H)}=\frac{y-1}{x-1}\prod_{i=1}^{4}(x-1)^{r(G_{n})-r(H_{i})}(y-1)^{n(H_i)}.
\end{equation*}

{\bf Subcase 3}. If $k(H)=\sum\limits_{i=1}^{4}k(H_{i})-5$, then we
can obtained, similarly, that
\begin{equation*}
(x-1)^{r(G_{n+1})-r(H)}(y-1)^{n(H)}=\frac{1}{(x-1)^2}\prod_{i=1}^{4}(x-1)^{r(G_{n})-r(H_{i})}(y-1)^{n(H_i)}.
\end{equation*}

{\bf Case 2}. $S=\emptyset$.

In this case, we consider the spanning subgraph $H$ of $G_{n+1}$,
which does not contain the new adding edge $e_n$, and
$|E(H)|=\sum\limits_{i=1}^{4}|E(H_i)|$.

{\bf Subcase 1}. If $k(H)=\sum\limits_{i=1}^{4}k(H_{i})-3$, then
\begin{equation*}
r(H) =
|V(H)|-k(H)=(4|V(G_{n})|-4)-(\sum_{i=1}^{4}k(H_{i})-3)=\sum_{i=1}^{4}r(H_{i})-1.
\end{equation*}
Moreover, we have
\begin{equation*}
n(H)=|E(H)|-r(H)=(\sum_{i=1}^{4}|E(H_i)|)-(\sum_{i=1}^{4}r(H_{i})-1)=\sum_{i=1}^{4}n(H_i)+1
\end{equation*}
and
\begin{equation*}
r(G_{n+1})-r(H)=(|V(G_{n+1})|-1)-(\sum_{i=1}^{4}r(H_i)-1)=\sum_{i=1}^{4}(r(G_{n})-r(H_i))
\end{equation*}
Thus,
\begin{equation*}
 (x-1)^{r(G_{n+1})-r(H)}(y-1)^{n(H)}=(y-1)\prod_{i=1}^{4}(x-1)^{r(G_{n})-r(H_{i})}(y-1)^{n(H_i)}.
\end{equation*}

{\bf Subcase 2}. If $k(H)=\sum\limits_{i=1}^{4}k(H_{i})-4$, then
\begin{equation*}
r(H)=|V(H)|-k(H)=(4|V(G_{n})|-4)-(\sum_{i=1}^{4}k(H_{i})-4)=\sum_{i=1}^{4}r(H_{i}).
\end{equation*}
Moreover, we have
$$n(H)=|E(H)|-r(H)=(\sum_{i=1}^{4}|E(H_i)|)-(\sum_{i=1}^{4}r(H_{i}))=\sum_{i=1}^{4}n(H_i)$$
and
$$r(G_{n+1})-r(H)=(|V(G_{n+1})|-1)-(\sum_{i=1}^{4}r(H_i))=\sum_{i=1}^{4}(r(G_{n})-r(H_i))-1$$
Thus we have
$$(x-1)^{r(G_{n+1})-r(H)}(y-1)^{n(H)}=\frac{1}{x-1}\prod_{i=1}^{4}(x-1)^{r(G_{n})-r(H_{i})}(y-1)^{n(H_i)}.$$

For the convenience of discussion, we use solid lines to join two
special vertices when the corresponding spanning subgraph of $G_n$
belongs to $\mathcal{G}_{1,n}$; Otherwise, we use dotted lines
instead of solid lines. We distinguish two types of the spanning
subgraph as shown in Figure 3.

\begin{table*}[!t]
\centering \caption{All combinations and corresponding
contributions.} \label{T1}
\begin{tabular}{cccc|cccc}
\toprule[1pt]
Configuration&S&Contribution&Type&Configuration&S&Contribution&Type\\
\midrule
\includegraphics[scale=0.36]{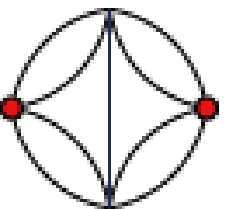}&$\{e_n\}$&$(y-1)^{2}T_{1,n}^{4}$&I&
\includegraphics[scale=0.36]{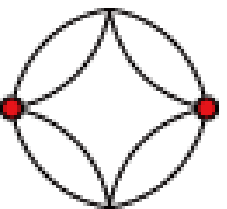}&$\emptyset$&$(y-1)T_{1,n}^{4}$&I\\
\includegraphics[scale=0.36]{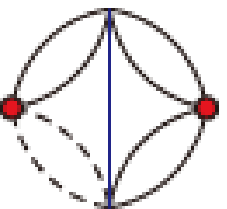}&$\{e_n\}$&$\frac{y-1}{x-1}T_{1,n}^{3}T_{2,n}$&I&
\includegraphics[scale=0.36]{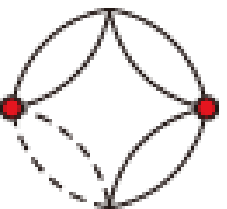}&$\emptyset$&$\frac{1}{x-1}T_{1,n}^{3}T_{2,n}$&I\\
\includegraphics[scale=0.36]{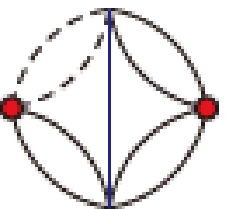}&$\{e_n\}$&$\frac{y-1}{x-1}T_{1,n}^{3}T_{2,n}$&I&
\includegraphics[scale=0.36]{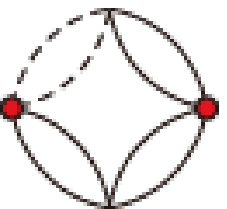}&$\emptyset$&$\frac{1}{x-1}T_{1,n}^{3}T_{2,n}$&I\\
\includegraphics[scale=0.36]{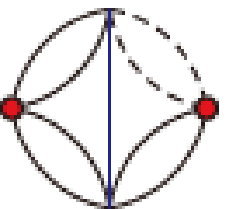}&$\{e_n\}$&$\frac{y-1}{x-1}T_{1,n}^{3}T_{2,n}$&I&
\includegraphics[scale=0.36]{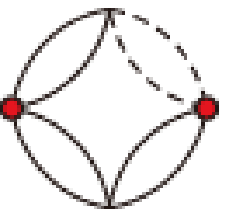}&$\emptyset$&$\frac{1}{x-1}T_{1,n}^{3}T_{2,n}$&I\\
\includegraphics[scale=0.36]{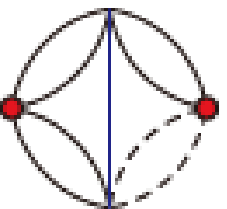}&$\{e_n\}$&$\frac{y-1}{x-1}T_{1,n}^{3}T_{2,n}$&I&
\includegraphics[scale=0.36]{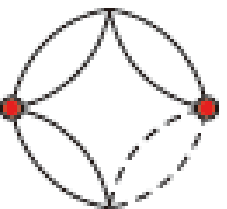}&$\emptyset$&$\frac{1}{x-1}T_{1,n}^{3}T_{2,n}$&I\\
\includegraphics[scale=0.36]{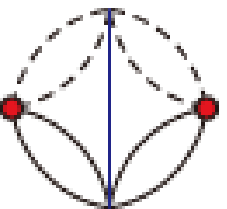}&$\{e_n\}$&$\frac{1}{(x-1)^{2}}T_{1,n}^{2}T_{2,n}^{2}$&I&
\includegraphics[scale=0.36]{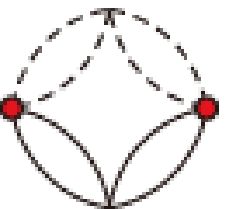}&$\emptyset$&$\frac{1}{x-1}T_{1,n}^{2}T_{2,n}^{2}$&I\\
\includegraphics[scale=0.36]{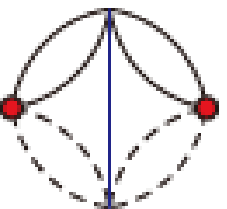}&$\{e_n\}$&$\frac{1}{(x-1)^{2}}T_{1,n}^{2}T_{2,n}^{2}$&I&
\includegraphics[scale=0.36]{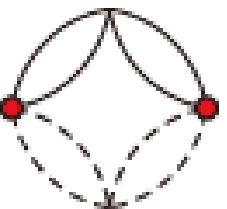}&$\emptyset$&$\frac{1}{x-1}T_{1,n}^{2}T_{2,n}^{2}$&I\\
\includegraphics[scale=0.36]{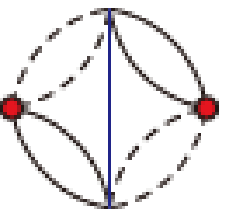}&$\{e_n\}$&$\frac{1}{(x-1)^{2}}T_{1,n}^{2}T_{2,n}^{2}$&I&
\includegraphics[scale=0.36]{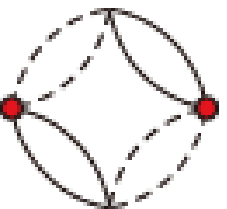}&$\emptyset$&$\frac{1}{x-1}T_{1,n}^{2}T_{2,n}^{2}$&II\\
\includegraphics[scale=0.36]{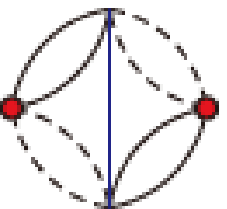}&$\{e_n\}$&$\frac{1}{(x-1)^{2}}T_{1,n}^{2}T_{2,n}^{2}$&I&
\includegraphics[scale=0.36]{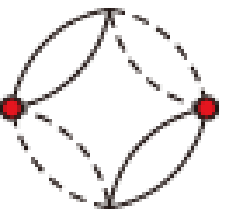}&$\emptyset$&$\frac{1}{x-1}T_{1,n}^{2}T_{2,n}^{2}$&II\\
\includegraphics[scale=0.36]{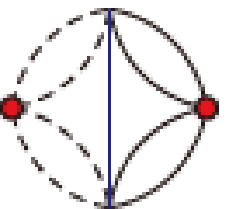}&$\{e_n\}$&$\frac{y-1}{x-1}T_{1,n}^{2}T_{2,n}^{2}$&II&
\includegraphics[scale=0.36]{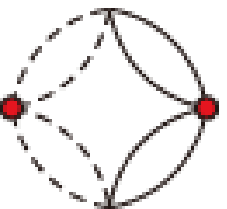}&$\emptyset$&$\frac{1}{x-1}T_{1,n}^{2}T_{2,n}^{2}$&II\\
\includegraphics[scale=0.36]{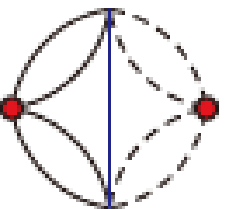}&$\{e_n\}$&$\frac{y-1}{x-1}T_{1,n}^{2}T_{2,n}^{2}$&II&
\includegraphics[scale=0.36]{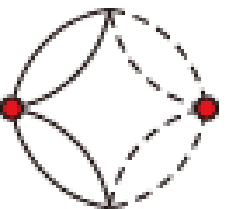}&$\emptyset$&$\frac{1}{x-1}T_{1,n}^{2}T_{2,n}^{2}$&II\\
\includegraphics[scale=0.36]{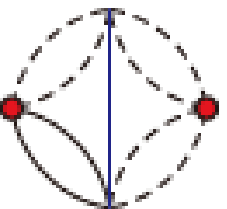}&$\{e_n\}$&$\frac{1}{(x-1)^{2}}T_{1,n}T_{2,n}^{3}$&II&
\includegraphics[scale=0.36]{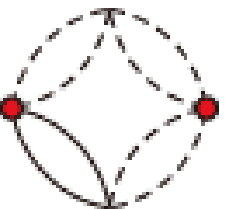}&$\emptyset$&$\frac{1}{x-1}T_{1,n}T_{2,n}^{3}$&II\\
\includegraphics[scale=0.36]{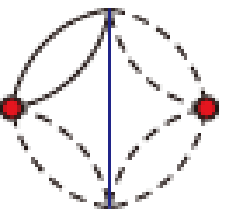}&$\{e_n\}$&$\frac{1}{(x-1)^{2}}T_{1,n}T_{2,n}^{3}$&II&
\includegraphics[scale=0.36]{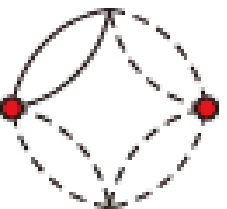}&$\emptyset$&$\frac{1}{x-1}T_{1,n}T_{2,n}^{3}$&II\\
\includegraphics[scale=0.36]{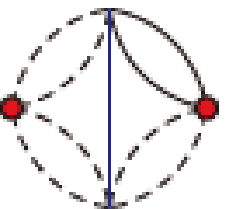}&$\{e_n\}$&$\frac{1}{(x-1)^{2}}T_{1,n}T_{2,n}^{3}$&II&
\includegraphics[scale=0.36]{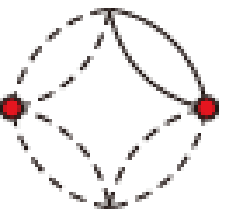}&$\emptyset$&$\frac{1}{x-1}T_{1,n}T_{2,n}^{3}$&II\\
\includegraphics[scale=0.36]{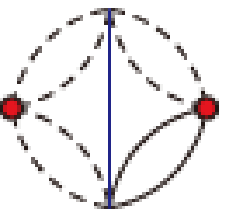}&$\{e_n\}$&$\frac{1}{(x-1)^{2}}T_{1,n}T_{2,n}^{3}$&II&
\includegraphics[scale=0.36]{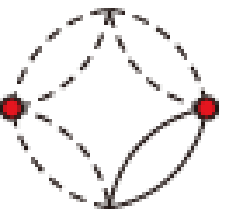}&$\emptyset$&$\frac{1}{x-1}T_{1,n}T_{2,n}^{3}$&II\\
\includegraphics[scale=0.36]{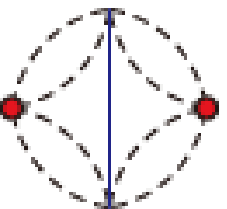}&$\{e_n\}$&$\frac{1}{(x-1)^{2}}T_{2,n}^{4}$&II&
\includegraphics[scale=0.36]{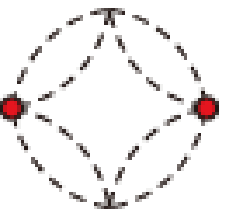}&$\emptyset$&$\frac{1}{x-1}T_{2,n}^{4}$&II\\
\bottomrule[1pt]
\end{tabular}
\end{table*}
%%%%%%%%%%%%%%%%%%%%%%%%%%%%%%%%%%%%%%%%%%%%

\begin{theorem}\label{t1}
The Tutte polynomial $T_{n+1}(x,y)$ of $G_{n+1}$ is given by
\begin{equation*}
T_{n+1}(x,y)=T_{1,n+1}(x,y)+T_{2,n+1}(x,y),
\end{equation*}
where the polynomials $T_{1,n+1}(x,y)$ and $T_{2,n+1}(x,y)$ satisfy
the following recursive relations:
\begin{eqnarray}\label{eq-1}
T_{1,n+1}(x,y)=y(y-1)T_{1,n}^4
+\frac{4y}{x-1}T_{1,n}^{3}T_{2,n}+\frac{2x+2}{(x-1)^{2}}T_{1,n}^{2}T_{2,n}^2,
\end{eqnarray}
\begin{eqnarray}\label{eq-2}
T_{2,n+1}(x,y)=\frac{2y+2}{{x-1}}T_{1,n}^{2}T_{2,n}^2+\frac{4x}{(x-1)^{2}}T_{1,n}T_{2,n}^{3}+\frac{x}{(x-1)^{2}}T_{2,n}^4
\end{eqnarray}
with initial conditions
$
T_{1,0}(x,y)=1,\ T_{2,0}(x,y)=x-1.
$
\end{theorem}
{\bf Proof}. The initial conditions are easily verified. The
strategy of the proof is to study all possible configurations of the
spanning subgraph $H_i$ in $G_n^i$ ($i=1,2,3,4$), and analyze which
kind of contributions they give to $T_{1,n}(x,y)$ and
$T_{2,n}(x,y)$. As shown in Table 1, a configuration produces a
basic term of form $T_{1,n}^{i}T_{2,n}^{j}(i+j=4)$, and by the
previous analysis, each basic term has to be multiplied by a factor
$(y-1)^2$, $\frac{y-1}{x-1}$, $\frac{1}{(y-1)^2}$ or $y-1$,
$\frac{1}{x-1}$ according to Case 1 and Case 2. More precisely, (i)
the factor is $(y-1)^2$ if $S=\{e_n\}$ and
$k(H)=\sum\limits_{i=1}^{4}k(H_{i})-3$; (ii) the factor is
$\frac{y-1}{x-1}$ if $S=\{e_n\}$ and
$k(H)=\sum\limits_{i=1}^{4}k(H_{i})-4$; (iii) the factor is
$\frac{1}{(x-1)^2}$ if $S=\{e_n\}$ and
$k(H)=\sum\limits_{i=1}^{4}k(H_{i})-5$; (iv) the factor is $(y-1)$
if $S=\emptyset$ and $k(H)=\sum\limits_{i=1}^{4}k(H_{i})-3$; (v) the
factor is $\frac{1}{x-1}$ if $S=\emptyset$ and
$k(H)=\sum\limits_{i=1}^{4}k(H_{i})-4$.

From Table 1, we can establish Eqs.(\ref{eq-1})-(\ref{eq-2}), and
the proof is completed. \hfill $\Box$ \\

According to Eq.(\ref{eq-2}), it is easy to prove by induction that
$x-1$ divides $T_{2,n}(x,y)$ in $\mathbb{Z}[x,y]$. Thus, we can
rewrite $T_{2,n}(x,y)$ as $(x-1)N_{n}(x,y)$ in $\mathbb{Z}[x,y]$,
and Theorem \ref{t1} can be reduced to the following:

\begin{theorem}\label{t2}
The Tutte polynomial $T_{n+1}(x,y)$ of $G_{n+1}$ is given by
\begin{equation*}
T_{n+1}(x,y)=T_{1,n+1}(x,y)+(x-1)N_{n+1}(x,y),
\end{equation*}
where the polynomial $T_{1,n+1}(x,y)$, $N_{n+1}(x,y)$ satisfy the
following recursive relations:
\begin{eqnarray*}
T_{1,n+1}(x,y) = y(y-1)T_{1,n}^{4}+4yT_{1,n}^{3}N_{n}+(2x+2)T_{1,n}^{2}N_{n}^{2},\\
N_{n+1}(x,y) =
(2y+2)T_{1,n}^{2}N_{n}^{2}+4xT_{1,n}N_{n}^{3}+x(x-1)N_{n}^{4}
\end{eqnarray*}
with initial conditions
$
T_{1,0}(x,y)=1,\ N_{0}(x,y)=1.
$
\end{theorem}

It is well-known that the evaluation of the Tutte polynomial for a
particular for a particular point at $(X,Y)$-plane is related to
some combinatorial information and algebraic properties of the graph
considered.

(1) $T(G;1,0)=$ the number of acyclic root-connected trees of $G$;

(2) $T(G;1,1)=\tau(G)$, i.e., the number of spanning trees;

(3) $T(G;0,1)=$ the number of indegree sequences of strongly
connected orientations of $G$;

(4) $T(G;-1,-1)=(-1)^{|E(G)|}(-2)^{dim(\mathcal{B})}$, where
$\mathcal{B}$ is the bicyclic space of $G$.

\begin{theorem}\label{t3}
The number of acyclic root-connected orientations of $G_{n}$ is
given by $T_{n}(1,0)=\prod\limits_{i=0}^{n}(i+1)^{2\times4^{n-i}};$
The number of indegree sequences of strongly connected orientations
of $G_{n}$ is given by
$T_{n}(0,1)=\frac{n}{2}\prod\limits_{i=0}^{n}(i+1)^{2\times4^{n-i}}$.
\end{theorem}
{\bf Proof}. By taking $x=1$ and $y=0$ in Theorem \ref{t2}, we have
$T_{n}(1,0)=T_{1,n}(1,0)$, and
\begin{equation}\label{eq-3}
T_{1,n}(1,0)=4T_{1,n-1}^{2}(1,0)N_{n}^{2}(1,0),
\end{equation}
\begin{equation}\label{eq-4}
N_{n}(1,0)=2T_{1,n-1}^{2}(1,0)N_{n-1}^{2}(1,0)+4T_{1,n-1}(1,0)N_{n-1}^{3}(1,0)
\end{equation}
A useful relation yields from Eqs.(\ref{eq-3}) and (\ref{eq-4})
\begin{equation}\label{eq-5}
\frac{N_{n}(1,0)}{T_{1,n}(1,0)}=\frac{1}{2}+\frac{N_{n-1}(1,0)}{T_{1,n-1}(1,0)}.
\end{equation}
It implies that
\begin{equation*}
\frac{N_{n}(1,0)}{T_{1,n}(1,0)}=\frac{n}{2}+\frac{N_{0}(1,0)}{T_{1,0}(1,0)},
\end{equation*}
and
\begin{equation}\label{eq-6}
N_{n}(1,0)=\frac{n+2}{2}T_{1,n}(1,0)
\end{equation}
since $T_{0}(1,0)=1$, $N_{0}(1,0)=1$.

Substituting (\ref{eq-6}) into (\ref{eq-3}) and using the initial
condition $T_{1,0}(1,0)=1$, we have
\begin{equation*}
T_{1,n}(1,0)=(n+1)^{2}T_{1,n-1}^{4}(1,0)=\prod_{i=0}^{n}(i+1)^{2\times4^{n-i}}.
\end{equation*}

Similarly, by taking $x=0$ and $y=1$ in Theorem \ref{t2}, we have
$T_{n}(0,1)=T_{1,n}(0,1)-N_{n}(0,1)$, and
\begin{eqnarray*}
T_{1,n}(0,1) &=& 4T_{1,n-1}^{3}(0,1)N_{n}(0,1)+2T_{1,n-1}^{2}(0,1)N_{n-1}^{2}(0,1),\\
N_{n}(0,1) &=& 4T_{1,n-1}^{2}(0,1)N_{n-1}^{2}(0,1).
\end{eqnarray*}
Using the same techniques, we can obtained
\begin{equation*}
T_{1,n}(0,1)=\frac{n+2}{2}N_{n}(0,1),\,\,\,
N_{n}(0,1)=\prod_{i=0}^{n}(i+1)^{2\times4^{n-i}}.
\end{equation*}
And
\begin{equation*}
T_{n}(0,1)=\frac{n+2}{2}N_{n}(0,1)-N_{n}(0,1)=\frac{n}{2}N_{n}(0,1)=\frac{n}{2}\prod_{i=0}^{n}(i+1)^{2\times4^{n-i}}.
\end{equation*} \hfill $\Box$ \\

From above, the number $T_{n}(0,1)$ of acyclic root-connected trees
and the number $T_{n}(1,0)$ of indegree sequences of strongly
connected orientations in $G_n$ are satisfied $T_{n}(0,1)
=\frac{n}{2}T_{n}(1,0)$.

\begin{theorem}\label{t4}
For a positive integer $n\geqslant 1$, the Tutte polynomial of
$G_n$, $T_n(x,y)$ along the line $y=x$ is given by
\begin{equation*}
T_n(x,x)=x(x^2+5x+2)^{\frac{4^n-1}{3}}.
\end{equation*}
\end{theorem}
{\bf Proof}. By taking $y=x$ in Theorem \ref{t2}, we have
\begin{equation*}
T_{1,n}(x,x)=x(x-1)T_{1,n-1}^4
+4xT_{1,n-1}^{3}N_{n-1}+(2x+2)T_{1,n-1}^{2}N_{n-1}^{2},
\end{equation*}
\begin{equation}\label{eq-7}
N_{n}(x,x)=(2x+2)T_{1,n-1}^{2}N_{n-1}^{2}+4xT_{1,n-1}N_{n-1}^{3}+x(x-1)N_{n-1}^{4},
\end{equation}
and $T_{1,0}(x,x)=N_{0}(x,x)=1$. It can be obtained easily that
$T_{1,n}(x,x)=N_{n}(x,x)$ by induction. Substituting it into
Eq.(\ref{eq-7}) and using the initial condition $N_{0}(x,x)=1$, we
have
\begin{equation*}
N_{n}(x,x)=(x^2+5x+2)N_{n-1}^{4}=(x^2+5x+2)^{\frac{4^n-1}{3}}.
\end{equation*}
Thus,
$T_{n}(x,x)=T_{1,n}(x,x)+(x-1)N_{n}(x,x)=xN_{n}(x,x)=x(x^2+5x+2)^{\frac{4^n-1}{3}}$. \hfill $\Box$ \\

Since the number of spanning trees is $\tau(G)=T(G;1,1)$, from
Theorem \ref{t4}, we can obtain immediately the following result,
which was also obtained in \cite{Zhang11PRE} by employing the
decimation teachnique.

\begin{corollary}\label{c5}
The number of spanning trees of $G_n$ is given by
$T_{n}(1,1)=8^{\frac{4^n-1}{3}}=2^{4^n-1}$. The asymptotic growth
constant is $$
\lim_{n\rightarrow\infty}\frac{ln\tau(G_n)}{|V(G_n)|}=\frac{3}{2}ln{2}\simeq1.0397.
$$
\end{corollary}

Similarly, $T_n(-1,-1)=(-1)\times(-2)^{\frac{4^n-1}{3}}=
(-1)^{\frac{4^{n+1}-1}{3}}(-2)^{\frac{4^n-1}{3}}=(-1)^{|E(G_n)|}(-2)^{dim(\mathcal{B})}$
by taking $x=-1$ in Theorem \ref{t4}. So, we have

\begin{corollary}
The dimension of the bicycle space of $G_n$ is $\frac{4^{n}-1}{3}$.
\end{corollary}

\appendixpage
\section{Tutte polynomial of two other scale-free networks}

In this section, we consider the Tutte polynomial of  two scale-free
networks (2,2)-flower and (1,3)-flower which have some similar
characteristics: identical degree sequences, without crossing edges
and always connected. Most of the topological properties of
(2,2)-flower and (1,3)-flower can be determined exactly
\cite{zzzcg}.

\subsection{(2,2)-flower}
If the newly added edge is ignored at each iterative generation,
then the fractal lattice $G_n$ considered above become the
(2,2)-flower $F_n$, see Figure 4. And the contributions to
$T_{1,n}(x,y)$ and $T_{2,n}(x,y)$ are degraded into the case of
$S=\emptyset$, and listed on the right of Table 1. So, we can obtain

%%%%%%%%%%%%%%%%%%%%%%%%%%%%%%%%%%%%%%%%%%%%%%%%%%%%%%%%%%%%%%%%%%55
                          %Figure 4$
%%%%%%%%%%%%%%%%%%%%%%%%%%%%%%%%%%%%%%%%%%%%%%%%%%%%%%%%%%%%%%%%%%%%55
\begin{figure}[ht!]
\begin{center}
\includegraphics[width=12cm]{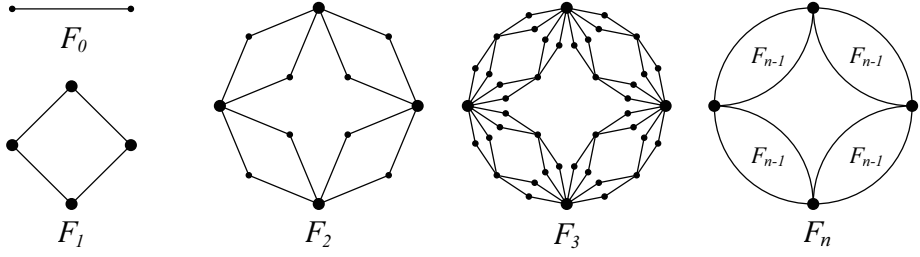}
\caption{Construction method of the (2,2)-flower $F_{n}$.}
\end{center}
\end{figure}
%%%%%%%%%%%%%%%%%%%%%%%%%%%%%%%%%%%%%%%%%%%%%%%%%%%%%%%%%%%%%%%%%%%%55
%%%%%%%%%%%%%%%%%%%%%%%%%%%%%%%%%%%%%%%%%%%%%%%%%%%%%%%%%%%%%%%%%%%55

\begin{equation*}
T_{n}(x,y)=T_{1,n}(x,y)+(x-1)N_{n}(x,y),
\end{equation*}
where the polynomials
$T_{1,n}(x,y)$, $N_n(x,y)$ satisfy the following recursive
relations:
\begin{eqnarray*}
T_{1,n}(x,y) &=& (y-1)T_{1,n-1}^4 +4T_{1,n-1}^{3}N_{n-1}+2(x-1)T_{1,n-1}^{2}N_{n-1}^{2},\\
N_{n}(x,y)   &=& 4T_{1,n-1}^{2}N_{n-1}^{2}+4(x-1)T_{1,n-1}N_{n-1}^{3}+(x-1)^{2}N_{n-1}^{4}.
\end{eqnarray*}

If $x=y=1$, then $T_n(1,1)=T_{1,n}(1,1)$ and
$T_{1,n}(1,1)=N_{n}(1,1)$. Thus, $T_n(1,1)=4T_{n-1}^{4}(1,1)$. Since
the initial value $T_{0}(1,1)=1$, we can obtain
$$\tau(F_{n})=T_{n}(1,1)=2^{\frac{2}{3}(4^n-1)}$$
and
$$
\lim_{n\rightarrow\infty}\frac{ln\tau(F_n)}{|V(F_n)|}=ln{2}\simeq0.6931.
$$

\subsection{(1,3)-flower}
%%%%%%%%%%%%%%%%%%%%%%%%%%%%%%%%%%%%%%%%%%%%%%%%%%%%%%%%%%%%%%%%%%55
                          %Figure 5$
%%%%%%%%%%%%%%%%%%%%%%%%%%%%%%%%%%%%%%%%%%%%%%%%%%%%%%%%%%%%%%%%%%%%55
\begin{figure}[ht!]
\begin{center}
\includegraphics[width=12cm]{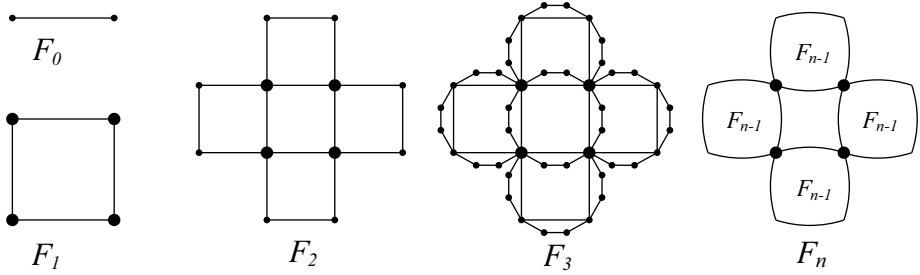}
\caption{Construction method of the (1,3)-flower $F_{n}$.}
\end{center}
\end{figure}
%%%%%%%%%%%%%%%%%%%%%%%%%%%%%%%%%%%%%%%%%%%%%%%%%%%%%%%%%%%%%%%%%%%%55
%%%%%%%%%%%%%%%%%%%%%%%%%%%%%%%%%%%%%%%%%%%%%%%%%%%%%%%%%%%%%%%%%%%55

For the (1,3)-flower $F_n$, see Figure 5, choosing two adjacent vertices
 with largest degree
 as special vertices in each iteration, we can obtain similarly
the following recursive relation:
$$
T_{n}(x,y)=T_{1,n}(x,y)+T_{2,n}(x,y),
$$
where the polynomials $T_{1,n}(x,y)$ and $T_{2,n}(x,y)$ satisfy the
following recursive relations:
\begin{eqnarray*}
T_{1,n}(x,y)&=& (y-1)T_{1,n-1}^4 +\frac{4}{x-1}T_{1,n-1}^{3}T_{2,n-1}+\frac{3}{x-1}T_{1,n-1}^{2}T_{2,n-1}^2+\frac{1}{x-1}T_{1,n-1}T_{2,n-1}^{3}, \\
T_{2,n}(x,y)&=& \frac{3}{{x-1}}T_{1,n-1}^{2}T_{2,n-1}^2+\frac{3}{x-1}T_{1,n-1}T_{2,n-1}^{3}+\frac{1}{x-1}T_{2,n-1}^4.
\end{eqnarray*}
with the initial conditions $ T_{1,0}(x,y)=1$, $T_{2,0}(x,y)=x-1$.
And
$$
T_{n}(x,y)=T_{1,n}(x,y)+(x-1)N_{n}(x,y),
$$
where the polynomials $T_{1,n}(x,y)$ and $N_n(x,y)$ satisfy the
following recursive relations:
\begin{eqnarray*}
T_{1,n}(x,y) &=& (y-1)T_{1,n-1}^4 +4T_{1,n-1}^{3}N_{n-1}+3(x-1)T_{1,n-1}^{2}N_{n-1}^{2}+(x-1)^{2}T_{1,n-1}N_{n-1}^{3},\\
N_{n}(x,y)   &=&
3T_{1,n-1}^{2}N_{n-1}^{2}+3(x-1)T_{1,n-1}N_{n-1}^{3}+(x-1)^{2}N_{n-1}^{4}
\end{eqnarray*}
with the initial conditions $T_{1,0}(x,y)=1$ and $N_0(x,y)=1$.

Similarly, if $x=y=1$, then $ T_n(1,1)=T_{1,n}(1,1)$ and
\begin{equation}\label{eq-8}
T_{1,n}(1,1)=4T_{1,n-1}^{3}(1,1)N_{n-1}(1,1),
\end{equation}
\begin{equation}\label{eq-9}
N_{n}(1,1)=3T_{1,n-1}^{2}(1,1)N_{n-1}^{2}(1,1).
\end{equation}
From Eqs.(\ref{eq-8}) and (\ref{eq-9}), we have
$$T_{1,n}(1,1)=(\frac{4}{3})^{n}N_{n}(1,1).$$
Since the initial value $N_{0}(1,1)=1$, we can obtain
$$\tau(F_{n})=T_{n}(1,1)=3^{(4^n-3n-1)/9}4^{(2\times4^n+3n-2)/9}$$
and
$$
\lim_{n\rightarrow\infty}\frac{ln\tau(F_n)}{|V(F_n)|}=\frac{1}{6}(4ln{2}+ln{3})\simeq0.6452.
$$
which are coincides with the results in \cite{lwzc} based on the
relationship between determinants of submatrices of the Laplacian
matrix.

\end{document}